%% file: linescircles.tex
\documentclass{gtart}
\usepackage{amsmath,amsfonts,amssymb,epsf,color}
\usepackage[xdvi]{graphicx}
\newcommand{\R}{\mathbb R}
\newcommand{\C}{\mathbb C}
\newcommand{\Z}{\mathbb Z}
\newcommand{\Q}{\mathbb Q}

\newcommand{\Gw}{\omega}
\newcommand{\p}{\partial}
\newcommand{\lk}{\operatorname{lk}}
\newcommand{\Tors}{\operatorname{Tors}}
\newcommand{\sign}{\operatorname{sign}}
\newcommand{\sminus}{\smallsetminus}

\theoremstyle{plain}

\newtheorem{thm}{Theorem}

\newtheorem{Lemma}{Lemma}

\newtheorem{ToR}{Theorem 1 Reformulated}

\newtheorem{CorOO}{Corollary 1 of Theorem 1}
 
\newtheorem{CorTO}{Corollary 2 of Theorem 1}
  
\newtheorem{CorOT}{Corollary 1 of Theorem 2}
 
\newtheorem{CorTT}{Corollary 2 of Theorem 2}
  
\newtheorem{Gth}{High-Dimensional Theorem}

\theoremstyle{definition}

\title{Lines and circles joining components of a link}
\author{Julia Viro (Drobotukhina)}
\address{Department of Mathematics, Uppsala University, S-751 06
Uppsala, Sweden}
\email{julia@math.uu.se}
\begin{document}

\begin{abstract}
We estimate from below the number of lines meeting each of given 4
disjoint smooth closed curves in a given cyclic order in the real 
projective 3-space and in a given linear order in $\R^3$.
Similarly, we estimate the number of circles meeting in a given cyclic
order given 6 disjoint smooth closed curves in Euclidean 3-space. 
The estimations are formulated in terms of linking numbers of the
curves and obtained by orienting of the corresponding
configuration spaces and evaluating of their signatures. 
This involves a study of a surface swept by
lines meeting 3 given disjoint smooth closed curves and a
surface swept in the 3-space by circles meeting 5 given disjoint smooth
closed curves. Higher dimensional generalizations of these results are
outlined.
\end{abstract}

\maketitle


\section{Introduction}\label{s1}

Consider several disjoint smooth closed curves in the 3-space.
If they are not linked to each other, one can pull them apart by an 
isotopy in such a way that any line joining two of 
them would not meet any other one. Of course, any three of the curves 
can be still connected by a circle, since a circle can be drawn through 
any three non-colinear points, but one can eliminate circles connecting
more than three of the curves.

However, as it is proved below, if the curves are linked well enough 
to each other, any 4 of them can be connected by a line, and any 6, 
by a circle. Moreover,  the numbers of such lines and circles are 
estimated from below in terms of pairwise linking numbers of the curves. 

\subsection{In Projective 3--Space}\label{s1.1}

{\bfseries Linking Numbers.\/}
Recall that linking number $\lk(A,B)$ is defined for any pair of 
disjoint oriented closed curves $A$, $B$ 
in an oriented closed 3--manifold $M$ realizing homology classes $[A],
[B]\in H_1(M)$ of finite order. It is defined as 
$\frac1n B\circ C$, where $C$ is a smooth chain transversal to $B$ 
with boundary $\partial C=nA$, and $\circ$ denotes the intersection 
number. Modulo one, linking number depends only on the homology classes of
the curves. This defines a well-known symmetric bilinear form 
$\Tors H_1(M)\times \Tors H_1(M)\to \Q/\Z$, which is a part of
Poincar\'e duality.

In particular, linking numbers of curves in the three--dimensional real
projective space $\R P^3$ are integers or
half-integers. Half-integer appears iff both curves are not
zero-homologous in $\R P^3$.

\begin{thm}[On Projective Lines Meeting 4 Curves, Pedestrian Version]\label{T1} 
Let $C_1$, $C_2$, $C_3$, $C_4$ be disjoint oriented smooth closed curves 
in $\R P^3$ 
in general position. Each real projective line $L$ intersecting the curves
in this  cyclic order can be equipped with a weight $\omega(L)=\pm 1$
such that 
$$\sum_{L}\omega(L)=
2\left(\lk(C_1,C_2)\lk(C_3,C_4) -\lk(C_2,C_3)\lk(C_4,C_1)\right),$$
where the sum runs over the set of all projective lines $L$ intersecting 
curves $C_1$, $C_2$, $C_3$, $C_4$ in this cyclic order.
\end{thm}

{\bfseries General Position Conditions and Their Roles.\/} In fact 
there are three general position conditions on $C_1$, $C_2$, $C_3$, 
$C_4$ which are important: \begin{enumerate} 
\item The number of lines intersecting all of them is finite. 
\item There is no line which meets the union $\cup_{i=1}^4C_i$ in five 
points.
\item For each line $L$ meeting each of $C_i$, the lines $L_i$ 
tangent  to $C_i$ at point $L\cap C_i$ are pairwise non-coplanar.  
\item For each line $L$ meeting each of $C_i$, there is no surface of
degree 2 containing $L$ and tangent to $C_i$ at $p_i$ for each
$i=1,2,3,4$.  
\end{enumerate}
The first condition is crucial for making the sum of weights finite.
The second, third and fourth conditions are not that necessary, but without 
accepting them we would not be able to claim that the weight of each
line is $\pm 1$. 

{\bf Set of Projective Lines Visiting Sets.} Let $A$, $B$, $C$ and $D$ be
subsets of a real projective space. Denote by $\mathcal P(A,B,C,D)$ the
set of projective lines which pass through points $a\in A$, $b\in B$,
$c\in C$, $d\in D$ in this cyclic order. 

Recall that a real projective line is homeomorphic to circle. 
A set of four points on a circle inherits two {\it cyclic\/} orders, 
each of which corresponds to an orientation of the circle. 
The set $\mathcal P(A,B,C,D)$ consists of those projective line 
for which one of these two cyclic orders induced from the line 
on $\{a,b,c,d\}$ coincides with the prescribed cyclic order.

{\bf Interpretation of Theorem \ref{T1} From Viewpoint of Differential 
Topology.} The
general position assumptions imply that $\mathcal P(C_1,C_2,C_3,C_4)$
is a 0--dimensional manifold. Theorem \ref{T1} claims that this manifold
has a natural orientation. Recall that an orientation of 0--manifold is
nothing but a function on it with values $\pm 1$. The sum of all the
values of an orientation $\omega$ of a compact 0--manifold $M$ is the
signature $\sign(M,\omega)$ of $M$. Hence Theorem \ref{T1} admits the 
following reformulation.

\begin{ToR}\label{ToR}
For any disjoint oriented smooth closed curves $C_1, C_2, C_3,
C_4\subset\R P^3$  
in general position, the set $\mathcal P=\mathcal P(C_1,C_2,C_3,C_4)$ is
a 0--manifold, which 
admits a natural orientation $\omega$ with signature 
$$\sign(\mathcal P,\omega)= 
2\left(\lk(C_1,C_2)\lk(C_3,C_4) -\lk(C_2,C_3)\lk(C_4,C_1)\right).$$ 
\end{ToR}

{\bfseries The Role of Order.\/}
In Theorem \ref{T1} curves $C_i$ are not assumed to be connected. 
$C_1\cup C_2\cup C_3\cup C_4$ can be considered as 
an oriented link in $\R P^3$ with components partitioned into four 
disjoint groups. A link with such partition is called a colored link. 
The colors (i.e., the groups of components) are assumed to be linearly 
ordered, and the order works in several ways. 

Obviously, only those projective lines are considered, which meet 
components of the link in the cyclic order of their colors.  
It does not make sense to speak about the linear order in which four 
points are positioned on a circle.   

However, the linear order of the curves cannot be replaced in Theorem
\ref{T1} by a cyclic order. The weights which the projective lines are
equipped with depend on the linear order. Even the expression for the 
sum of the weights changes sign under cyclic permutation $(1,2,3,4)$ of 
curves $C_i$.

{\bfseries Low Bounds for the Number of Lines Meeting Four
Curves.\/}
Theorem \ref{T1} gives low bounds for the number of lines intersecting 
four given closed oriented curves in  $\R P^3$ :

\begin{CorOO}\label{C1T1}
Under hypothesis of Theorem \ref{T1}, the number of real projective
lines intersecting curves $C_1$, $C_2$, $C_3$, $C_4$ in this  cyclic 
order (i.~e., the number of points in $\mathcal P(C_1,C_2,C_3,C_4)$) 
is at least
 $$2\left|\lk(C_1,C_2)\lk(C_3,C_4) -\lk(C_2,C_3)\lk(C_4,C_1)\right|.$$   
\end{CorOO}

Summing up estimates provided by this for all cyclic orders, we get the
following total estimate:

\begin{CorTO}\label{C2T1} 
Under hypothesis of Theorem \ref{T1}, the total number of real
projective lines intersecting each of $C_i$ is at least
$$4\max_{\substack{ (i,j,k,l)\text{ is a}\\
                    \text{permutation}\\
		    \text{of }(1,2,3,4)}} 
 \left|\lk(C_i,C_j)\lk(C_k,C_l) -\lk(C_j,C_k)\lk(C_l,C_i)\right|
$$
\end{CorTO}

{\bf The Role of Orientations of the Curves.}
Orientations of curves $C_i$ are needed in Theorem \ref{T1} for making 
linking numbers $\lk(C_i,C_j)$ defined. By changing the orientations, one
may hope to improve the low bound in some cases. However, if each of
$C_i$ is connected, the low bound does not depend on the orientations.
Indeed, reversing the orientation of any single curve changes the signs
of both summands in the expression $\lk(C_1,C_2)\lk(C_3,C_4) -
\lk(C_2,C_3)\lk(C_4,C_1)$.

\begin{figure}[thb]
\centerline{\input{4lines+.pstex_t}\quad \input{4lines-.pstex_t}\quad
\input{4lines0.pstex_t}\ }
\caption{Links consisting of four disjoint lines in $\R P^3$,
which represent all three isotopy classes of such non-oriented links.
The links shown on the left hand side and in the middle are mirror
images of each other. The link on the right hand side is amphicheiral. 
Orientations are shown to make $\lk(C_i,C_j)$ defined.}
\label{4lines}
\end{figure}
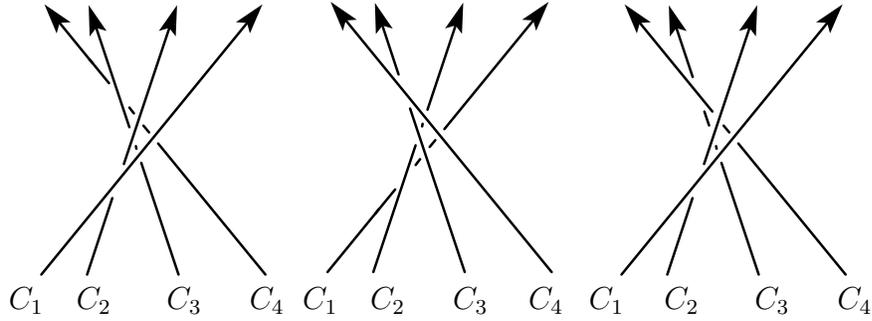

{\bf Examples: Links Made of 4 Lines.}
There are three isotopy classes of links in $\R P^3$ made of four
disjoint projective lines, see \cite{DV}. Their representatives are
shown in Figure \ref{4lines}.

In the link shown on the left hand side, $\lk(C_i,C_j)=+\frac12$ for any
$i$ and $j$, in the link shown in the middle $\lk(C_i,C_j)=-\frac12$ for
any $i$ and $j$. Hence, in both cases Corollary 2 of Theorem 1 gives
the trivial low bound for the number of lines meeting each of $C_i$.
{\it This low bound is exact: in both cases the isotopy class
contains a link for which there is no line meeting all for components.}

\begin{figure}[bht]
\centerline{\input{4linesOnH+.pstex_t}}
\caption{Two links made of 4 lines on a hyperboloid.}
\label{4linesOnH}
\end{figure}
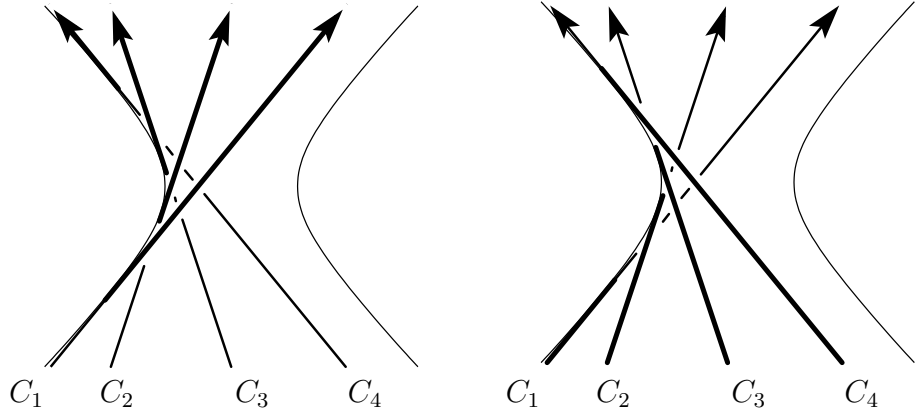

\begin{figure}[bht]
\centerline{\input{lineOffH.pstex_t}}
\caption{Moving line $L$ off hyperboloid $H\supset L$. In the projective
space such a move cannot be made in a single plane. Here it is presented as a 
composition of two rotations in different planes.
Draw a plane $P$ containing $L$, see (b). Plane $P$ is tangent to $H$ 
at some point. Rotate line $L$ in plane $P$ around the point, see (c). 
The result $M$ is a line tangent to $H$. Draw another plane $Q$ 
containing $M$ and intersecting $H$ in a non-singular conic, see (d). 
Rotate $M$ in plane $Q$ around a point $p\in M\sminus H$. }
\label{lineOffH}
\end{figure}
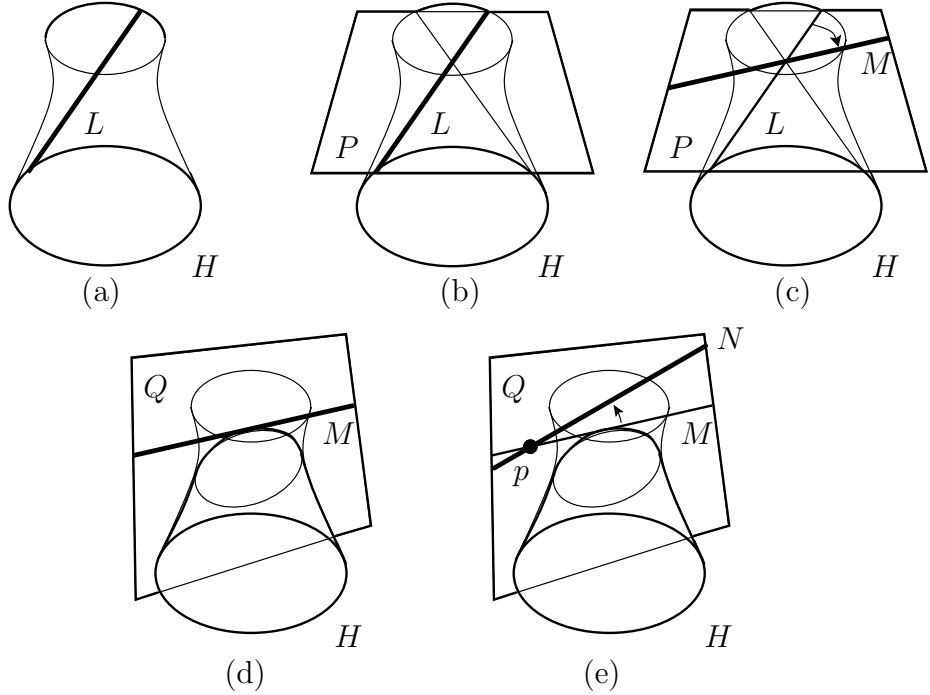

Indeed, both configuration of lines can be placed on a hyperboloid, as a
collection of 4 lines belonging to a family of pairwise disjoint lines
which cover the whole hyperboloid. See Figure \ref{4linesOnH}. For this
representative of the isotopy class there are infinitely many lines
meeting all 4 components. Namely, each line of the other family of lines
covering the hyperboloid meets each line of the first  family of lines.
However, one can move one component of the link into one of the domains
bounded by the hyperboloid. See Figure \ref{lineOffH}. Any line which
meets each of the other three lines must lie on the hyperboloid, because
it meets it in at least three points. Hence it has no chance to meet the
line which was taken off the hyperboloid.

In the rightmost link of Figure \ref{4lines}, $\lk(C_3,C_4)=-\frac12$
and $\lk(C_i,C_j)=+\frac12$ for any other values of $i$, $j$. Hence,
Corollary 2 of Theorem 1 claims that there are at least
$4|\frac12(-\frac12)-\frac12\frac12|=4|-2\frac14|=2$ lines meeting lines
$C_1$, $C_2$, $C_3$ and $C_4$. According to Corollary 1, at least one of
them meet the lines in cyclic order $C_1$, $C_2$, $C_3$, $C_4$, and at
least one, in cyclic order $C_1$, $C_2$, $C_4$, $C_3$. 

\begin{figure}[htb]
\centerline{\input{4lines0andH.pstex_t}}
\caption{Amphicheiral link of 4 lines and the hyperboloid swept by lines
meeting three of its components.}
\label{4lines0andH}
\end{figure}
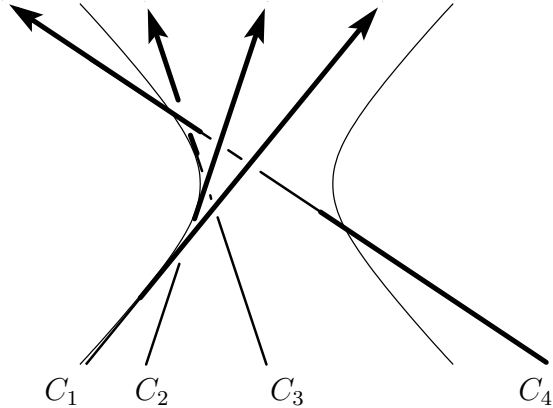

{\it These bounds are also exact.\/} 
To see this, consider the hyperboloid $H$ swept by lines meeting 
$C_1$, $C_2$ and $C_3$. See Figure \ref{4lines0andH}.
The fourth line, $C_4$, intersects the hyperboloid in two points. The
intersection points are in two connected components of
$H\sminus(C_1\cup C_2\cup C_3)$ which are adjacent to $C_3$. This is seen
in Figure \ref{4lines0andH}, and this is the only possibility. Indeed, 
if $C_4$ did not meet $H$ at all, by an isotopy inverse to the one shown
in Figure \ref{lineOffH} one could put entirely on $H$, but then the
linking numbers would be as in one of two other cases shown in Figure
\ref{4linesOnH}. If the intersection points were in the same connected
component of $H\sminus(C_1\cup C_2\cup C_3)$, moving the points towards
each other one would construct an isotopy of the link with $C_4$ tangent
to $H$. This would give again one of two other links. Finally, one of
the points of $C_4\cap H$ cannot be in the component of 
$H\sminus(C_1\cup C_2\cup C_3)$ bounded by $C_1$ and $C_2$, because then 
the equality $\lk(C_1,C_4)=\lk(C_2,C_4)$ would not hold. In fact, 
these arguments prove the isotopy classification of links in $\R P^3$
consisting of four lines. For details see \cite{DV}. 

The lines meeting $C_1$, $C_2$, $C_3$, and $C_4$ lie on $H$ and pass 
through the intersection points of $C_4$ and $H$. 
The line passing through the point which lies on $H$ between 
$C_2$ and $C_3$ meets the lines in cyclic order $C_1$, $C_2$, $C_4$, $C_3$. 
The line passing through the point which lies on $H$ between 
$C_3$ and $C_1$ meets the lines in cyclic order $C_1$, $C_2$, $C_3$,
$C_4$.

{\bf Absence of Upper Bounds.}
Under hypothesis of Theorem 1, there is no upper bound for the number of
real projective lines intersecting curves $C_1$, $C_2$, $C_3$ and $C_4$ 
in terms of linking numbers or any other invariants, which do not change 
under isotopy of the curves. 

Indeed, one can draw a line $L$ connecting a point on $C_1$ with a point 
on $C_2$, then move $C_3$ and $C_4$ by an isotopy to make them meeting $L$, 
so that $L$ would meet $C_1$, $C_2$ and $C_3$ in a prescribed cyclic order. 
Then one can make short pieces of $C_1$, $C_2$ and $C_3$ near their 
intersection points with $L$ line segments. The set of lines meeting 
three given lines covers a quadric surface. A strip of this surface is 
formed by lines meeting
the rectilinear pieces of $C_1$, $C_2$ and $C_3$. 
By a small isotopy of $C_4$ 
in a neighborhood of $C_4\cap L$ one can force $C_4$ to meet  
this strip arbitrary number of times near 
the point  $C_4\cap L$. Each of the intersection points is on
one of the lines which form the strip. These lines meet $C_1$, $C_2$,
$C_3$ and $C_4$ in the prescribed cyclic order.

Upper bounds exist in the case when the curves are defined by algebraic
equations and admit complexification. Then the complexification provides
terms in which upper bounds can be formulated. We will return to this,
but we need, first, to consider the corresponding complex problems.

{\bf Similar Problem from Complex Algebraic Geometry.}
What happens to the problems discussed above, if we pass from smooth 
curves in $\R P^3$ to complex curves in the complex projective 3--space 
$\C P^3$? For generic $C_1$, $C_2$, $C_3$ and $C_4$ there is 
finitely many complex projective lines meeting each of $C_i$. Points on 
a complex projective line do not have any distinguished cyclic order. Hence,
it does not make sense to consider lines which meet the curves in
whatever cyclic order. However one may wonder about the total number of
all complex projective lines intersecting each of
four given disjoint generic algebraic curves $C_1$, $C_2$, $C_3$, $C_4$.
This is a classical problem of complex enumerative algebraic geometry,
and its solution is well-known.

As this is the complex algebraic geometry, estimates are replaced by
an exact value. It is equal to $2d_1d_2d_3d_4$, where $d_i$ is the 
order of $C_i$, that is the intersection number of $C_i$ with a generic 
hyperplane. (If the curves are not in general position, there may be
infinitely many lines, or the number of lines counted with positive
multiplicities is still $2d_1d_2d_3d_4$.)

To prove this, remark, first, that the number of lines meeting 
each of $C_1$, $C_2$, $C_3$, and $C_4$ equals the number of 
intersection points of $C_4$ and the surface swept by lines meeting 
each of curves $C_1$, $C_2$ and $C_3$. 
By the Bezout theorem this number is the product of $d_4$ by the degree
of the surface. Therefore it depends linearly on $d_4$. By symmetry, 
it depends linearly on each of $d_i$. Thus, this is $kd_1d_2d_3d_4$. 
In the case, when each of the curves $C_i$ is a line,
the surface swept by lines meeting each of $C_1$, $C_2$, $C_3$ is a
hyperboloid, hence $k=2$.\qed 

{\bf In Real Algebraic Geometry.}
If under assumptions of Theorem \ref{T1} curves $C_1$, $C_2$, $C_3$ and
$C_4$ are algebraic, by taking complexification of all the varieties
involved we obtain the complex algebraic situation considered above.
Some of $2d_1d_2d_3d_4$ complex projective lines meeting the
complexifications of $C_1$, $C_2$, $C_3$ and $C_4$ may be imaginary (not
real), but the complexification of any real line meeting  $C_1$, $C_2$, 
$C_3$ and $C_4$ is among these $2d_1d_2d_3d_4$ complex projective lines. 
Thus, in contrast to purely topological situation discussed above, in
real algebraic situation there is an upper bound:  
the total number of real lines meeting real algebraic curves $C_1$,
$C_2$, $C_3$ and $C_4$ in $\R P^3$ is at
most $2d_1d_2d_3d_4$, where $d_i$ is the order of $C_i$.  

\subsection{In Affine 3-Space}\label{s1.3}
Results on links in $\R P^3$ are applied to links in $\R^3$, just
because $\R^3$ is embedded in $\R P^3$. Besides, there appear new 
opportunities related to natural linear orders of points on a line.

Let $A$, $B$, $C$ and $D$ be subsets of an affine space. Denote by 
$\mathcal A(A,B,C,D)$ the
set of affine lines which pass through points $a\in A$, $b\in B$,
$c\in C$, $d\in D$ in this linear order. 

Notice that sets $A$, $B$, $C$, $D$ are not assumed to be pairwise
distinct. 
 
\begin{thm}[On Affine Lines Meeting 4 Curves]\label{T2}
For any disjoint oriented smooth closed curves $C_1, C_2, C_3,
C_4\subset\R^3$  
in general position, the set $\mathcal A=\mathcal A(C_1,C_2,C_3,C_4)$ 
of lines meeting curves $C_1$, $C_2$, $C_3$, $C_4$ in their order
is a 0--manifold, which admits a natural orientation $\omega$ with signature 
$$\sign(\mathcal A,\omega)= \lk(C_1,C_2)\lk(C_3,C_4).$$ 
\end{thm}

As I was informed by Oleg Viro, this theorem was proven using
different arguments by Michael Polyak (unpublished).

Almost everything written in the preceding section about Theorem \ref{T1} 
can be repeated here concerning Theorem \ref{T2}.
In particular, the general position conditions are exactly the same.

Notice that in Theorem \ref{T2} lines under consideration meet the
curves in a fixed {\it linear\/} order. 
There are two linear orders for points on a line corresponding to
orientations of the line, and here any
of these two orders is meant. The transition from cyclic to linear 
order is responsible for the difference between expressions for the
sum of weights in theorems \ref{T1} and \ref{T2}. Theorem \ref{T2}
applied to all linear orders corresponding to the same cyclic order
implies Theorem \ref{T1} for the case when all the curves are contained
in an affine part of the projective space.

{\bf Low Bounds for the Number of Lines Meeting Four Curves.} As the
absolute value of signature of a 0--manifold cannot exceed the number of
points, Theorem \ref{T2} gives low bounds for the number of lines 
intersecting four given closed oriented curves in $\R^3$:

\begin{CorOT}\label{C1T2}
Under hypothesis of Theorem \ref{T2}, the number of lines intersecting
curves  $C_1$, $C_2$, $C_3$, $C_4$ in this order is at least 
 $$|\lk(C_1, C_2)\lk(C_3,C_4)|.$$ 
\end{CorOT}

Summing up the estimates provided by this for all orders, we get the
following total estimate:

\begin{CorTT}\label{C2T2}
Under hypothesis of Theorem \ref{T2}, the total number of lines intersecting
curves  $C_1$, $C_2$, $C_3$, $C_4$ is at least 
 $$4\left(|\lk(C_1, C_2)\lk(C_3,C_4)|+|\lk(C_1, C_3)\lk(C_2,C_4)|+
 |\lk(C_1,C_4)\lk(C_2,C_3)|\right).$$ 
\end{CorTT}

{\bf Connecting Two or Three Curves.} The next two theorems are specific
for affine case, their projective counter-parts are trivial.

\begin{thm}[On Affine Lines Meeting 2 Curves Intermittingly]\label{T3}
For any disjoint oriented smooth closed curves  $C_1$ and $C_2$
in $\R^3$ in general position, the set $\mathcal A=
\mathcal A(C_1,C_2,C_1,C_2)$
is a 0--dimensional manifold, which admits an orientation $\omega$ with
signature
$$\sign(\mathcal A,\omega)=(\lk(C_1,C_2))^2.$$
\end{thm}

\begin{thm}[On Affine Lines Meeting 3 Curves Intermittingly]\label{T4}
For any disjoint oriented smooth closed curves  $C_1$, $C_2$ and $C_3$
in $\R^3$ in general position, the sets $\mathcal A_{1213}=\mathcal
A(C_1,C_2,C_1,C_3)$ and $\mathcal A_{1231}=\mathcal A(C_1,C_2,C_3,C_1)$
are 0--manifolds, which
admit orientations $\omega_1$ and $\omega_2$, respectively, such that
$$\sign(\mathcal A_{1213},\omega_1)=\sign(\mathcal A_{1231},\omega_2)
=\lk(C_1,C_2)\lk(C_1,C_3).$$
\end{thm}

\subsection{Connecting by Circles}\label{s1.4}
In theorems which follow the role that is played by projective or affine 
lines above passes to circles. The space of all circles in $\R^3$ is a
manifold of dimension 6. A condition that a circle meets a fixed generic
curve specifies a hypersurface in the space of all circles. A
transversal intersection of 6 such hypersurfaces is a 0-dimensional
space. This is why the number of curves is increased to 6.

Let $A_1$, \dots, $A_n$ be subsets of an Euclidean space. Denote by 
$\mathcal S(A_1,\dots,A_n)$ the set of circles and lines which pass
through points $a_1\in A_1$, \dots, $a_n\in A_n$ in this cyclic order.

\begin{thm}[On Circles Meeting 6 Curves]\label{T5} 
For any disjoint oriented smooth closed  curves  $C_1$, $C_2$, $C_3$, 
$C_4$, $C_5$, $C_6$ in $\R^3$ in general position, the set 
$\mathcal S=\mathcal S(C_1, C_2, C_3, C_4, C_5, C_6)$ is a 0--manifold,
which admits an orientation $\omega$ with signature
$$\sign(\mathcal S,\omega)=
\lk(C_1,C_2)\lk(C_3,C_4)\lk(C_5,C_6) -\lk(C_2,C_3)\lk(C_4,C_5)
\lk(C_6,C_1).$$
\end{thm}                                                           

\begin{thm}[Circles Meeting 5 Curves]\label{T6}
For any disjoint oriented smooth closed  curves  $C_1$, $C_2$, $C_3$, 
$C_4$, $C_5$ in $\R^3$ in general position, sets 
\begin{align*} 
\mathcal S_{121345}&=\mathcal S(C_1, C_2, C_1, C_3, C_4, C_5)\\
\mathcal S_{123145}&=\mathcal S(C_1, C_2, C_3, C_1, C_4, C_5) 
\end{align*}
are 0--manifolds,
which admit orientations $\omega_1$ and $\omega_2$, respectively, such
that
\begin{align*}
\sign(\mathcal S_{121345},&\omega_1)= \\
&\lk(C_1,C_2)\lk(C_1,C_3)\lk(C_4,C_5) -\lk(C_2,C_1)\lk(C_3,C_4)
\lk(C_5,C_1)\\
\sign(\mathcal S_{123145},&\omega_2)= \\
&\lk(C_1,C_2)\lk(C_3,C_1)\lk(C_4,C_5) -\lk(C_2,C_3)\lk(C_1,C_4)
\lk(C_5,C_1).
\end{align*}  
\end{thm}

\begin{thm}[Circles Meeting 4 Curves]\label{T7}
For any disjoint oriented smooth closed  curves  $C_1$, $C_2$, $C_3$, 
$C_4$ in $\R^3$ in general position, sets 
\begin{align*} 
\mathcal S_{121234}&=\mathcal S(C_1, C_2, C_1, C_1, C_3, C_4)\\
\mathcal S_{121324}&=\mathcal S(C_1, C_2, C_1, C_3, C_2, C_4)\\ 
\mathcal S_{123124}&=\mathcal S(C_1, C_2, C_3, C_1, C_2, C_4)\\ 
\end{align*}
are 0--manifolds,
which admit orientations $\omega_1$, $\omega_2$ and $\omega_3$, respectively, 
such that
\begin{align*}
\sign(\mathcal S_{121234},&\omega_1)= \\
&(\lk(C_1,C_2))^2\lk(C_3,C_4) -\lk(C_2,C_1)\lk(C_2,C_3)
\lk(C_4,C_1)\\
\sign(\mathcal S_{121324},&\omega_2)= \\
&\lk(C_1,C_2)\lk(C_1,C_3)\lk(C_2,C_4) -\lk(C_2,C_1)\lk(C_3,C_2)
\lk(C_4,C_1)\\
\sign(\mathcal S_{123124},&\omega_3)= \\ 
&\lk(C_1,C_2)\lk(C_3,C_1)\lk(C_2,C_4) -\lk(C_2,C_3)\lk(C_1,C_2)
\lk(C_4,C_1).
\end{align*}  
\end{thm}

\subsection{In High-Dimensional Spaces}\label{s1.5}
The results presented above can be generalized to submanifolds of 
projective and Euclidean spaces of higher dimensions. Namely,

\begin{Gth}\label{GT1}
For any disjoint oriented smooth closed submanifolds $C_1$, $C_2$, $C_3$,
$C_4$ of $\R P^{2n+1}$ in general position with $\dim C_1=\dim C_3=p$, 
$\dim C_2=\dim C_4=q$ and $p+q=2n$, the set 
$\mathcal P=\mathcal P(C_1,C_2,C_3,C_4)$ is a 0--manifold, which 
admits a natural orientation $\omega$ with signature 
$$\sign(\mathcal P,\omega)= 
2\left(\lk(C_1,C_2)\lk(C_3,C_4) -\lk(C_2,C_3)\lk(C_4,C_1)\right).$$ 
\end{Gth}

Odd dimension of the projective space is necessary for its orientability
and existence of linking numbers.

\begin{Gth}\label{GT2}
For any disjoint oriented smooth closed submanifolds $C_1$, $C_2$, $C_3$,
$C_4$ of $\R^n$ in general position with $\dim C_1+\dim C_2=n-1$ and
$\dim C_3+\dim C_4=n-1$, 
the set $\mathcal A=\mathcal A(C_1,C_2,C_3,C_4)$ 
of lines meeting curves $C_1$, $C_2$, $C_3$, $C_4$ in their order
is a 0--manifold, which admits a natural orientation $\omega$ with signature 
$$\sign(\mathcal A,\omega)= \lk(C_1,C_2)\lk(C_3,C_4).$$ 
\end{Gth}

\begin{Gth}\label{GT3}
For any disjoint oriented smooth closed submanifolds $C_1$ and $C_2$
of $\R^n$ in general position with $\dim C_1+\dim C_2=n-1$, the set 
$\mathcal A=\mathcal A(C_1,C_2,C_1,C_2)$
is a 0--dimensional manifold, which admits an orientation $\omega$ with
signature
$$\sign(\mathcal A,\omega)=(\lk(C_1,C_2))^2.$$  
\end{Gth}

\begin{Gth}\label{GT4}
For any disjoint oriented smooth closed submanifolds $C_1$, $C_2$ and $C_3$
of $\R^n$ in general position with $\dim C_2=\dim C_3=n-1-\dim C_1$, 
the sets $\mathcal A_{1213}=\mathcal
A(C_1,C_2,C_1,C_3)$ and $\mathcal A_{1231}=\mathcal A(C_1,C_2,C_3,C_1)$
are 0--manifolds, which
admit orientations $\omega_1$ and $\omega_2$, respectively, such that
$$\sign(\mathcal A_{1213},\omega_1)=\sign(\mathcal A_{1231},\omega_2)
=\lk(C_1,C_2)\lk(C_1,C_3).$$
\end{Gth}

\begin{Gth}\label{GT5}
For any disjoint oriented smooth closed submanifolds $C_1$, $C_2$, $C_3$, 
$C_4$, $C_5$, $C_6$ of $\R^n$ in general position with 
$\dim C_1=\dim C_3=\dim C_5=p$,
$\dim C_2=\dim C_4=\dim C_6=q$ and $p+q=n-1$, the set 
$\mathcal S=\mathcal S(C_1, C_2, C_3, C_4, C_5, C_6)$ is a 0--manifold,
which admits an orientation $\omega$ with signature
$$\sign(\mathcal S,\omega)=
\lk(C_1,C_2)\lk(C_3,C_4)\lk(C_5,C_6) -\lk(C_2,C_3)\lk(C_4,C_5)
\lk(C_6,C_1).$$
\end{Gth}

\begin{Gth}\label{GT6}
For any disjoint oriented smooth closed submanifolds $C_1$, $C_2$, $C_3$, 
$C_4$, $C_5$ of $\R^n$ in general position with $\dim C_1=\dim C_4=p$,
$\dim C_2=\dim C_3=\dim C_5=q$ and $p+q=n-1$, the set
$$ \mathcal S_{121345}=\mathcal S(C_1, C_2, C_1, C_3, C_4, C_5)$$
is a 0--manifold with natural orientation $\omega$ such that
\begin{align*} 
 \sign(\mathcal S_{121345},&\omega)= \\
&\lk(C_1,C_2)\lk(C_1,C_3)\lk(C_4,C_5) -\lk(C_2,C_1)\lk(C_3,C_4)
\lk(C_5,C_1).
\end{align*}
If, moreover, $p=q$, the set   
$$
\mathcal S_{123145}=\mathcal S(C_1, C_2, C_3, C_1, C_4, C_5) 
$$
is a 0--manifold,
which has a natural orientation $\omega_1$  such
that
\begin{align*}
\sign(\mathcal S_{123145},&\omega_1)= \\
&\lk(C_1,C_2)\lk(C_3,C_1)\lk(C_4,C_5) -\lk(C_2,C_3)\lk(C_1,C_4)
\lk(C_5,C_1).
\end{align*}  
\end{Gth}

\begin{Gth}\label{GT7}
For any disjoint oriented smooth closed submanifolds $C_1$, $C_2$, $C_3$, 
$C_4$ in $\R^n$ in general position with $\dim C_1=\dim C_3=p$, $\dim
C_2=\dim C_4=q$ and $p+q=n-1$, the set 
$$
\mathcal S_{121234}=\mathcal S(C_1, C_2, C_1, C_2, C_3, C_4)$$
 is a 0--manifold with natural orientation $\omega_1$ such that 
$$\sign(\mathcal S_{121234},\omega)= 
(\lk(C_1,C_2))^2\lk(C_3,C_4) -\lk(C_2,C_1)\lk(C_2,C_3)\lk(C_4,C_1).$$
 If, moreover, $p=q$, the sets    
\begin{align*}
\mathcal S_{121324}&=\mathcal S(C_1, C_2, C_1, C_3, C_2, C_4)\\ 
\mathcal S_{123124}&=\mathcal S(C_1, C_2, C_3, C_1, C_2, C_4)\\ 
\end{align*}
are 0--manifolds,
which admit orientations $\omega_1$ and $\omega_2$,  respectively, 
such that
\begin{align*}
\sign(\mathcal S_{121324},&\omega_1)= \\
&\lk(C_1,C_2)\lk(C_1,C_3)\lk(C_2,C_4) -\lk(C_2,C_1)\lk(C_3,C_2)
\lk(C_4,C_1)\\
\sign(\mathcal S_{123124},&\omega_2)= \\ 
&\lk(C_1,C_2)\lk(C_3,C_1)\lk(C_2,C_4) -\lk(C_2,C_3)\lk(C_1,C_2)
\lk(C_4,C_1).
\end{align*}  
\end{Gth}

\subsection{Quadrisecants in Literature}\label{s1.6}

I am not aware about any publication, where the same problems were 
considered. 
However, there are several papers devoted to similar problems about
quadrisecants of knots and links in $\R^3$. A {\it quadrisecant\/} of a link
$L\subset \R^3$ is a line meeting $L$ at least at four points. In other
words, a quadrisecant of a link $L$ is a point of $\mathcal A(L,L,L,L)$.  
In a sense, the problems about quadrisecants of a single curve are more 
sophisticated than the problems considered above. In particular, 
the problems about quadrisecants of a single curve do belong to low 
dimensional topology. 

The most classical of the papers on quadrisecants is a dissertation of 
Erika Pannwitz \cite{Pannwitz} published in 1933. Here are the main two 
results of that paper:

{\bf Pannwitz Theorem 1.} {\it
Any generic piecewise linear knot in $\R^3$ which is not isotopic to the 
unknot has a quadrisecant (i.e., a line meeting the knot at four points).}

 {\bf Pannwitz Theorem 2.} {\it 
Any generic piecewise linear link of two components $L_1\cup L_2\subset\R^3$ 
which are linked in the sense that each represents a non-trivial
homotopy class in the complement of the other has an essential quadrisecant, 
that is a line which intersects twice each of the components in such a way 
that the segment enclosed between the intersection points of the line with 
one component contains an intersection point with the other one.}

Besides, Pannwitz \cite{Pannwitz} proved low bounds for the number of 
quadrisecants in both cases.  

In 1982 H.R.Morton and D.M.Q.Mond \cite{MM} proved the same results in the 
differential category, that is for generic smooth knots and links 
(their statement about two-component links was weaker: 
they assumed that the components have nonvanishing linking number). 
This part of their results follows also from Theorem \ref{T3}.

In 1994 Greg Kuperberg \cite{Kuper} proved the following theorem
extening, in particular, these results to arbitrary tame knot and links.

{\bf Kuperberg Theorem.} {\it Every non-trivial tame link in $R^3$ has a
quadrisecant.}

A recent paper \cite{BCSDS} by Ryan Budney, James Conant, Kevin P. Scannell, 
and Dev Sinha contains a result on quadrisecants of a generic knot
that admits a formulation very close in its spirit to the results
of this paper.
Quadrisecants of a generic knot $K$ constitute a finite set. Some part of it
can be equipped with a natural orientation such that its signature is
the simplest non-trivaial Vassiliev-Goussarov invariant of $K$. This
invariant is of degree two and called the Conway knot invariant. 
 
\subsection{Possible extension of classical results to projective 
links}\label{s1.7}
Although Pannwitz-Morton-Mond Theorems formulated above 
make sense for links in the projective space $\R P^3$, 
both statements, and proofs require essential corrections. 

In the proofs, the notion of convex hull is used at a crucial point.
This notion is not applicable to a link in $\R P^3$.

\begin{figure}[ht]
\centerline{\includegraphics[scale=.8]{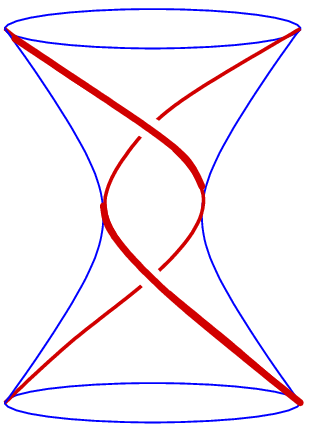}}
\caption{}
\label{f1}
\end{figure}

As for the statements, in the projective space there are two
isotopy types of knots which are obvious counter-parts for the unknot
in $\R^3$: the types represented by a non-singular non-empty conic
(e.g., circle) and projective line. Neither circle, nor slightly
perturbed projective line (say, one component non-singular plane cubic
curve) has a quadrisecant. However there exist knots which do not
belong to these isotopy classes and nonetheless have no quadrisecant.
For example, a knot $2_1$ from my table \cite{D2} of prime projective
knots with at most 6 crossings is isotopic to a real algebraic curve on
a hyperboloid of bidegree $(3,1)$.  In Figure \ref{f1} this curve is
shown, together with the hyperboloid.

A line which does not lie on the
hyperboloid meets the hyperboloid (and hence the curve) in at most two
points, while a line on the hyperboloid meets the curve at one or three
point (depending on which family of generatrices it belongs to). Thus
the curve has no quadrisecant. Of course, the mirror image of this 
curve has no quadrisecant, too. An easy modification of  arguments
used by  Pannwitz and Kuperberg, shows that: 

{\it Any tame knot in $\R P^3$ without quadrisecants belongs to one of 
the four isotopy types listed above. \/}

Even more changes are required in the case of links. There exist
infinitely many isotopy types of two-component links in the projective
space which have no essential quadrisecant. Indeed, this property has 
any link with components separated by a hyperboloid. 

\subsection{Gratitudes}\label{s1.8}
I am grateful to Oleg Viro for interest to my work and useful discussions.
I am also grateful  to Mathematics Department of Uppsala University and
MSRI for partial support of this work.

\section{Proofs}\label{s2}

Proofs of all theorems formulated above are based on the same ideas. In
the simplest and most profound way these ideas work in the proof of
Theorem \ref{T1}.
Therefore, we concentrate first on this theorem.  

\subsection{Surfaces of Secants}\label{s2.1}

Let $C_1$, $C_2$ and $C_3$ be three smooth disjoint closed curves 
in $\R P^3$. Similarly to the notation introduced above, denote by
$\mathcal P(C_1,C_2,C_3)$ the set of projective lines each of which 
passes through points $a\in C_1$, $b\in C_2$, $c\in C_3$.

Projective lines belonging to $\mathcal P(C_1,C_2,C_3)$ are called {\it
common secants \/} of $C_1$, $C_2$ and $C_3$.  The union of all common 
secants of   $C_1$, $C_2$ and $C_3$ is called the {\it secant surface\/}
of  $C_1$, $C_2$ and $C_3$ and denoted by $s(C_1,C_2,C_3)$. 
Thus 
$$s(C_1,C_2,C_3)=\bigcup_{L\in \mathcal P(C_1,C_2,C_3)}L.$$

The role of $s(C_1,C_2,C_3)$ becomes clear from the following simple remark: 
given any curve $C_4\subset \R P^3$, for each point
$x\in s(C_1,C_2,C_3)\cap C_4$ there exists a line passing through $x$
and meeting $C_1$, $C_2$, $C_3$, $C_4$. The cyclic order, in which this
line meets $C_1$, $C_2$, $C_3$, $C_4$, depends on position of $x$ on
$s(C_1,C_2,C_3)$ with respect to curves $C_1$, $C_2$ and $C_3$.
See Figure \ref{f2}.

\begin{figure}[htb]
\centerline{\input{3+1curves.pstex_t}}
\caption{}
\label{f2}
\end{figure}

{\bf Resolution of singularities.}
For most triples of curves $s(C_1,C_2,C_3)$ is not a 2-manifold. 
If the curves $C_i$ are algebraic, $s(C_1,C_2,C_3)$ is an algebraic surface, 
but usually with lots of singular points. In what follows, we use a 
standard trick to resolve most of the singularities. For generic curves, 
the trick resolves all singularities.

A common secant $L$ of curves $C_1$, $C_2$, and $C_3$ with 
three marked points $p_i$, $i=1,2,3$ chosen from $C_i\cap L$ is called a
{\it pointed secant.\/}
Denote by $S(C_1,C_2,C_3)$ the set 
$$\{(x,p_1,p_2,p_3)\in\R P^3\times C_1\times C_3\times
C_3\mid \text{ points }x,p_1,p_2,p_3 \text{ are colinear}\}$$
One may think on $S(C_1,C_2,C_3)$ as the union of
all pointed secants of $C_1$, $C_2$, and $C_3$. Let us call 
$S(C_1,C_2,C_3)$ the {\it pointed secant surface\/} of $C_1$, $C_2$,
$C_3$. Denote the {\it set\/} of all pointed secants by 
$T(C_1,C_2,C_3)$. Since a pointed secant is defined by the three marked
points, we can identify $T(C_1,C_2,C_3)$ with
$$\{(p_1,p_2,p_3)\in C_1\times C_2\times C_3 \mid  \text{ points }
p_1,p_2,p_3 \text{ are colinear} \} $$

Notice natural inclusions $$
\begin{aligned} 
s(C_1,C_2,C_3)&\subset\R P^3,\\  
 S(C_1,C_2,C_3)&\subset\left(\R P^3\right)^4,\\ 
T(C_1,C_2,C_3)&\subset \left(\R P^3\right)^3.
 \end{aligned}$$  
Using these inclusions, induce topology on $s(C_1,C_2,C_3)$, 
$S(C_1,C_2,C_3)$ and $T(C_1,C_2,C_3)$.
The natural map 
$$u:S(C_1,C_2,C_3)\to s(C_1,C_2,C_3) : (x,p_1,p_2,p_3)\mapsto x$$ 
forgetting the marked points is continuous with respect to 
these natural topological structures. 
Denote the composition of this map with the inclusion
$s(C_1,C_2,C_3)\hookrightarrow\R P^3$ by  $U$. Thus, $U$ is a
natural map $S(C_1,C_2,C_3)\to \R P^3$ acting by formula 
$(x,p_1,p_2,p_3)\mapsto x$ and having image $s(C_1,C_2,C_3)$. 

The natural map 
$$q:S(C_1,C_2,C_3)\to T(C_1,C_2,C_3) : (x,p_1,p_2,p_3)\mapsto (p_1,p_2,p_3)$$ 
is a locally trivial
fibration with fiber circle. The points marked on secants 
provide three disjoint sections of this fibration defined by
formulas $(p_1,p_2,p_3)\mapsto(p_i,p_1,p_2,p_3)$ with $i=1,2,3$.
The images of these
sections can be described by the following formulas
$$P_i=\{(x,p_1,p_2,p_3)\in S(C_1,C_2,C_3)\mid x=p_i\}.$$ 

A fiber $q^{-1}(p_1,p_2,p_3)$ of $q$ is naturally identified with the 
projective line in $\R P^3$ passing through (colinear) points 
$p_1$, $p_2$ and $p_3$. There exists a unique projective transformation 
$q^{-1}(p_1,p_2,p_3)\to\R P^1$ sending $p_1 \mapsto (1\,{:}\,0)$, $p_2\mapsto 
(1\,{:}\,1)$
and $p_3\mapsto (0\,{:}\,1)$. These transformations altogether define a
homeomorphism 
$$S(C_1,C_2,C_3)\to\R P^1\times T(C_1,C_2,C_3),$$
under which $P_1$, $P_2$ and $P_3$ are mapped to fibers 
$(1\,{:}\,0)\times T(C_1,C_2,C_3)$,
$(1\,{:}\,1)\times T(C_1,C_2,C_3)$, and  $(0\,{:}\,1)\times T(C_1,C_2,C_3)$,
respectively. 

\subsection{Genericity assumptions}\label{s2.2}
For a pointed secant  $(p_1,p_2,p_3)$ of $C_1$, $C_2$, $C_3$, 
denote the tangent line of $C_i$ at $p_i$ by $L_i$.
A pointed secant $(p_1,p_2,p_3)$ is said to be {\it regular,\/} 
if the lines $L_i$ with $i=1,2,3$ are pairwise disjoint. 

Along a regular pointed secant, a germ of $U(s(C_1,C_2,C_3))$ is
$C^1$-approximated by a germ of surface $s(L_1,L_2,L_3)$, which is a
hyperboloid. See Figure \ref{f3}.

\begin{figure}[htb]
\centerline{\input{reg_sec.pstex_t}}
\caption{}
\label{f3}
\end{figure}

If pointed secant $(p_1,p_2,p_3)$ is not regular, that is, say 
$L_i$, $L_j$ are coplanar for some 
$i,j\in\{1,2,3\}$ with $i\ne j$, then the projection of the curves $C_i$ 
and $C_j$ from the third of the marked points, say $p_k$, are tangent 
to each other at the images of $p_i$ and $p_j$. If this tangency is
quadratic and $L_k$ is coplanar neither to $L_i$, nor to $L_j$, then 
the pointed secant line $(p_1,p_2,p_3)$ is said to be {\it almost
regular\/} and $p_k$ is called its {\it special point.\/} See Figure
\ref{f4}.
\begin{figure}[b]
\centerline{\input{alm_reg_sec.pstex_t}}
\caption{}
\label{f4}
\end{figure}

The surface $s(L_1,L_2,L_3)$ is a union of two planes. It $C^1$-approximates 
the germ of $U$ along the pointed secant at all points of the secant except 
$p_k$. 
To get a sufficient jet at $p_k$, one can replace $C_k$ still with line 
$L_k$ tangent to $C_k$ at $p_k$, while the other two curves should be 
replaced by appropriate conics.            

A good model for a germ of $U(s(C_1,C_2,C_3))$ near an almost regular
pointed secant is Whitney umbrella. In fact Whitney umbrella is
$s(C_1,C_2,C_3)$, where $C_1$ is a line,  $C_2$ and $C_3$ are parabolas
with axes parallel to $C_1$ and symmetric to each other against $C_1$.
See Figure \ref{f5}. 

\begin{figure}
\centerline{\input{whitn_umbr.pstex_t}}
\caption{}
\label{f5}
\end{figure}

For arbitrary $C_1$, $C_2$ and $C_3$, along an almost regular pointed 
secant a germ of $U(s(C_1,C_2,C_3))$ is approximated by a germ of a
surface projectively equivalent to the Whitney umbrella.       

A triple $C_1$, $C_2$, $C_3$ of disjoint smooth closed 
curves in the projective space is called {\it regular,\/} if all but
finitely many of its pointed secants are regular and each non-regular
pointed secant is almost regular. 

\begin{Lemma}\label{L0}
Regularity is a generic property of a triple of disjoint smooth 
closed curves.
\end{Lemma}

For the notion of generic property, we refer to a general theory of
general position presented by Wall \cite{Wall}. Cf. \cite{Kuper}.
The proof is straightforward. In fact, Lemma \ref{L0} is almost 
entirely covered by Lemma 2.4 from the paper \cite{Kuper} of Kuperberg.
\qed                                    

To understand regularity and almost regularity of secants, take a
pointed secant $(p_1,p_2,p_3)$ and consider a view on $C_2$ and $C_3$
from viewpoint $p_1$. 

If this is a regular secant, along it we will see
crossing point of $C_2$ and $C_3$. If the same line is not involved in
other pointed secants, the crossing is simple transversal. 
It is stable: under a 
small move of $p_1$ along $C_1$ the secant moves, and the other two
points $p_2$ and $p_3$ move along $C_2$ and $C_3$ with speeds
determined by the speed of $p_1$.

If this is an almost regular secant and $p_1$ is its special point,
$C_2$ and $C_3$ will look tangent to each other quadratically. This
tangency is not stable: under a small move of $p_1$ in one direction
along $C_1$ the point of tangency splits into two crossings, under a move 
in the opposite direction it disappears. Passing $p_1$ along $C_1$ gives
rise to a second Reidemeister move.

{\bf Secants versus Reidemeister moves.}
This interpretation of regularity and almost regularity relates Lemma \ref{L0}
to facts well-known from the knot theory. Indeed, it is well-known that a
projection of two smooth curves from a generic point has only
transversal double points, and in a one parameter generic family only
Reidemeister moves  happen at isolated moments. Thus when we travel
along one of three generic curves under consideration, we meet only
regular secants at all but finite number of points. 
Since we are interested only in intersection points of images of two 
different branches, the first Reidemeister move of the picture seen from
our moving viewpoint is not of any value. 
The third Reidemeister moves correspond to the lines which are underlying 
for more than one pointed secant, each of which are regular. Second 
Reidemeister moves correspond to almost regular secants.

\begin{Lemma}[\bf Pointing resolves singularities]\label{L2}
For any regular triple $C_1$, $C_2$, $C_3$ of smooth closed
curves in $\R P^3$, the set $T(C_1,C_2,C_3)$ of pointed secants is a
smooth 1-dimensional submanifold of $\left(\R
P^3\right)^3$,
the pointed secant surface $S(C_1,C_2,C_3)$ is a
smooth two-dimensional submanifold of $\left(\R P^3\right)^4$ 
and $U:S(C_1,C_2,C_3)\to \R P^3$ is a differentiable map. 
The only singularities of $U$ are pinch points at the special 
points of non-regular pointed secants.   
\end{Lemma}

We skip the proof. It is a straightforward
application of Implicit Function Theorem and Whitney's 
characterization of a pinch point. \qed    

\begin{Lemma}[\bf Topology of pointed secant surface]\label{L3}
For any regular triple $C_1$, $C_2$, $C_3$ of smooth  closed
curves in $\R P^3$, each connected component of pointed secant 
surface $S(C_1,C_2,C_3)$ is diffeomorphic to torus $S^1\times S^1$
under diffeomorphism which maps each pointed secant line to a fiber
$pt\times S^1$, while the sections $P_1$, $P_2$, $P_3$  
of the fibration $S(C_1,C_2,C_3)\to T(C_1,C_2,C_3)$ corresponding
to the marked points
are mapped to three fibers of the complementary family, $S^1\times pt$. 
\end{Lemma} 

This is an immediate corollary of Lemma \ref{L2}. \qed 

\subsection{Orientation of pointed secant surface}\label{s2.3} 

{\bf Along a regular secant.} Let $C_1$, $C_2$ and $C_3$ be
oriented smooth closed curves in $\R P^3$ forming a regular
triple. Denote by $O_i$ the given orientation of $C_i$. (Orientation is
considered as a function taking values $\pm1$ on bases of each tangent 
space.) 

Let $(p_1,p_2,p_3)$ be a regular pointed secant line  of 
$C_1$, $C_2$ and $C_3$.  
Denote by $L$ the line in $\R P^3$ passing through $p_1$, $p_2$ and
$p_3$, and by $L_i$ the tangent line of $C_i$ at $p_i$.
Orient $L_i$ according to the orientation of $C_i$. Recall that
regularity of the pointed secant means that these lines are pairwise
disjoint. Therefore, each pair $L_i$, $L_j$ of them with $i\ne j$ 
has a well-defined linking number $\lk(L_i,L_j)$ equal to $\pm\frac12$.

Consider hyperboloid $s(L_1,L_2,L_3)$, which osculates secant surface
$$s(C_1,C_2,C_3)=U(S(C_1,C_2,C_3))$$ along $L$. For any
tangent vector $e$ of $S(C_1,C_2,C_3)$, its image $dU(e)$ is a tangent
vector of  $s(L_1,L_2,L_3)$. The line $L$ (as a line on a hyperboloid)
is a circle embedded two-sidedly in $s(L_1,L_2,L_3)$. Choose one of the
sides and choose a non-zero vector $e_i$ at $p_i$ tangent to $C_i$
and directed towards the chosen side of $L$ on $s(L_1,L_2,L_3)$.   

\begin{Lemma}\label{L4}
Let $\sigma$ be a permutation of $(1,2,3)$. Then 
$$\lk(L_{\sigma(2)},L_{\sigma(3)})O_{\sigma(1)}(e_{\sigma(1)})$$
does not depend on $\sigma$.
\end{Lemma}

\begin{proof} Denote line $L_i$ oriented along $e_i$ by
$\overline{L_i}$. 
Recall (see \cite{DV}) that the product
$\lk(L_1,L_2)\lk(L_2,L_3)\lk(L_3,L_1)$ of all pairwise linking numbers of 
three pairwise disjoint oriented lines in $\R P^3$ does not depend on
the orientations of the lines. In particular, 
 $$\lk(L_1,L_2)\lk(L_2,L_3)\lk(L_3,L_1)=
 \lk(\overline{L_1},\overline{L_2})\lk(\overline{L_2},\overline{L_3})
 \lk(\overline{L_3},\overline{L_1}).$$
All factors in the right hand side of this equality are equal, because
oriented links $\overline{L_i}\cup\overline{L_j}$ are isotopic to each
other. Hence
$$
4\lk(\overline{L_1},\overline{L_2})
\lk(\overline{L_2},\overline{L_3})\lk(\overline{L_3},\overline{L_1})=
 \lk(\overline{L_{\sigma(2)}},\overline{L_{\sigma(3)}}). 
$$
 Change of the orientation
of a line multiplies the linking number of this line with other line
by $-1$. Hence
$$\lk(\overline{L_{\sigma(2)}},\overline{L_{\sigma(3)}})= \lk(L_{\sigma(2)},L_{\sigma(3)})
O_{\sigma(2)}(e_{\sigma(2)})O_{\sigma(3)}(e_{\sigma(3)}).$$

Summarizing we get
$$4\lk(L_1,L_2)\lk(L_2,L_3)\lk(L_3,L_1)=
\lk(L_{\sigma(2)},L_{\sigma(3)})
O_{\sigma(2)}(e_{\sigma(2)})O_{\sigma(3)}(e_{\sigma(3)}).$$
Multiplying this equality by $O_{1}(e_{1})O_2(e_2)O_3(e_3)$ we obtain
an expression for 
$\lk(L_{\sigma(2)},L_{\sigma(3)})O_{\sigma(1)}(e_{\sigma(1)})$ 
independent of $\sigma$:
\begin{multline*}
\lk(L_{\sigma(2)},L_{\sigma(3)})O_{\sigma(1)}(e_{\sigma(1)})=\\
 4\lk(L_1,L_2)\lk(L_2,L_3)\lk(L_3,L_1)O_1(e_1)O_2(e_2)O_3(e_3).
\end{multline*}
\end{proof}

The given ordering of curves $C_i$, $i=1,2,3$ 
defines a cyclic order of points $p_i$, $i=1,2,3$, which, in turn,
defines an orientation of $L$. For each $i=1,2,3$, choose a non-zero 
vector $f_i$ at $p_i$ tangent to $L$ and directed along the orientation of
$L$.

Let $V$ be a connected neighborhood of $L$ in $S(C_1,C_2,C_3)$ disjoint
from non-regular pointed secants. Without loss of generality
we can assume that $V$ is a union of pointed secants. Since $V$ does not
meet non-regular pointed secants, the restriction of $U$ to $V$
is an immersion. Moreover, since the restriction of $U$ to $L$ is
embedding, we can choose $V$ small enough to make $U|_V$ an embedding.

Vectors $e_i,f_i$ form a basis of the tangent space of $V$ at $p_i$.
The bases $(e_1,f_1)$, $(e_2,f_2)$ and $(e_3,f_3)$ define the same 
orientation of $V$. Let us correct it by $\frac12\lk(L_2,L_3)O_1(e_1)$,
that is consider the orientation $O$ of $V$ which takes value 
$\frac12\lk(L_2,L_3)O_1(e_1)$ on $(e_1,f_1)$. By Lemma \ref{L4}, $O$
takes values $\frac12\lk(L_3,L_1)O_2(e_2)$ on $(e_2,f_2)$ and 
$\frac12\lk(L_1,L_2)O_3(e_3)$  on $(e_3,f_3)$. 

{\bf Extending the orientation across almost regular secants.}
This orientation can be defined at any point of a regular pointed 
secant, the construction depends continuously of the point.
Thus, we have defined an orientation of the complement of
non-regular pointed secants in $S(C_1,C_2,C_3)$. This orientation
depends on the orientation of each curve $C_i$ and a cyclic order of
these curves. 

\begin{Lemma}\label{L5}
Orientation $O$ extends across almost regular pointed secants to
the whole $S(C_1,C_2,C_3)$.
\end{Lemma}
                                                                    
\begin{proof} When $(p_1,p_2,p_3)$ moves in $T(C_1,C_2,C_3)$ 
and passes through an almost regular secant 
$(\bar p_1,\bar p_2,\bar p_3)$ with special point $\bar p_1$, 
both the linking number $\lk(L_2,L_3)$, and  $O_1(e_1)$ change and,
hence, their product does not change. 

Indeed, by the definition of almost regular pointed secant, 
the projections of $C_2$ and $C_3$ from $\bar p_1$ are tangent to each
other at the images of  $\bar p_2$ and $\bar p_3$. 
The tangency is quadratic, so, when the center of the projection 
$p_1$ moves along $C_1$ passing $\bar p_1$, the
projection of $C_2\cup C_3$ experiences the second Reidemeister move.
When $p_1$ is on one side of $\bar p_1$, there are 
two intersection points of the images of $C_2$ and $C_3$ which are close
to the image of non-regular secant under projection from $\bar p_1$, 
when $p_1$ gets to the other side, the intersection points disappear. 
The vanishing pair of intersection points
corresponds to a pair of pointed secant lines, say 
$(p_1,p_2^+,p_3^+)$ and $(p_1,p_2^-,p_3^-)$ 
close to $(\bar p_1,\bar p_2,\bar p_3)$. On curve $T(C_1,C_2,C_3)$ 
they are on the opposite sides of $(\bar p_1,\bar p_2,\bar p_3)$. 
Thus, $\bar p_1$ is
the turning point for the map $U|_{P_1}$, and $O_1(e_1)$ changes the
sign, when one jumps from $(p_1,p_2^+,p_3^+)$ to 
$(p_1,p_2^-,p_3^-)$. 

At the intersection points of the projections
of $C_2$ and $C_3$ from $p_1$ corresponding to these lines,
the writhe numbers are opposite to each other. These writhe number equals
the doubled linking number of the lines tangent to the branches.
Therefore the linking number changes when  $(p_1,p_2,p_3)$ 
passes $(\bar p_1,\bar p_2,\bar p_3)$.
\end{proof}

\subsection{Orientation of the pointed sections}\label{s2.4} As
above, let $C_1$, $C_2$ and $C_3$ be oriented smooth closed curves in 
$\R P^3$ forming a regular triple. 
Each fiber of the fibration $S(C_1,C_2,C_3)\to T(C_1,C_2,C_3)$ intersects 
sections $P_1$, $P_2$ and $P_3$ in three
points, $p_1$, $p_2$, $p_3$, respectively, and is divided by these 
intersection points into three arcs.
Denote by $\Sigma_{i,j}$ the union of closures of those arcs which 
connect $p_i$ with $p_j$. This is a compact surface with boundary  
$P_i\cup P_j$. 

Equip $\Sigma_{i,j}$ with the orientation
induced by $O$. Equip $P_1$, $P_2$ and $P_3$  with the orientations such
that $\p\Sigma_{1,2}=P_1\cup(-P_2)$, $\p\Sigma_{2,3}=P_2\cup(-P_3)$ and
$\p\Sigma_{3,1}=P_3\cup(-P_1)$.   

\subsection{Degree of map of a pointed section to curve}\label{s2.5}
\begin{Lemma}\label{L6}
The degree of map $U_i:P_i\to C_i$ defined by the natural map
$U:S(C_1,C_2,C_3)\to\R P^3$ is equal to $2\lk(C_j,C_k)$ with
$\{i,j,k\}=\{1,2,3\}$.
\end{Lemma}

\begin{proof} Consider the case $i=1$. Choose a point $p_1\in C_1$ such
that each common secant of $C_1$, $C_2$, $C_3$ passing through $p_1$ is
regular. This means that the projection of $C_2\cup C_3$ from $p_1$ is
generic and $p_1$ is a regular value of $U_1:P_1\to C_1$. The preimage
of $p_1$ under $U_1$ consists of points $(p_1,p_1,p_2^k,p_3^k)$ of
pointed secant lines $(p_1,p_2^k,p_3^k)$ passing through $p_1$.
So, we have to prove that the sum of local degrees of $U_1$ over all
these points is $2\lk(C_2,C_3)$. 

At point $(p_1,p_1,p_2^k,p_3^k)$ choose the vector $e_1$ such that basis
$(e_1,f_1)$ defines the orientation $O$. Since $f_1$ is outwards normal
vector for $\Sigma_{3,1}$ at $(p_1,p_1,p_2^k,p_3^k)$, vector $e_1$
defines the orientation of $P_1$ at this point. Hence the local degree
of $U_1$ at this point is $O_1(e_1)$.

Due to our choice of $e_1$, the value of
$O$ on $(e_1,f_1)$ is $+1$. On the other hand, this value, by the
definition of $O$, is $2\lk(L_2,L_3)O_1(e_1)$. Therefore 
$$O_1(e_1)=2\lk(L_2,L_3).$$ 
The left hand side is the local degree
of $U_1$ at $(p_1,p_1,p_2^k,p_3^k)$.
The right hand side is the local writhe of the
projection of $C_2\cup C_3$ from $p_1$ at the image of $p_2^k$.
The sum of local writhe numbers at all the intersection points of the
projections of $C_2$ and $C_3$ is equal to $2\lk(C_2,C_3)$, see, e.g.
\cite{D3}.
\end{proof}

\begin{Lemma}\label{L7}
The image of a fundamental cycle of $\Sigma_{3,1}$ under $U$ 
is a smooth singular
chain with boundary $2\left(\lk(C_1,C_2)[C_3]-\lk(C_2,C_3)[C_1]\right)$. 
\end{Lemma}
This is an immediate corollary of Lemma \ref{L6}.\qed

\begin{Lemma}\label{L8}
The intersection number of $U_*[\Sigma_{3,1}]$ with an oriented closed
curve $C_4$ disjoint from $C_1\cup C_2\cup C_3$ is
$$2\det 
\begin{pmatrix} 
\lk(C_1,C_2)&\lk(C_1,C_4)\\
\lk(C_3,C_2)&\lk(C_3,C_4)
\end{pmatrix}
$$
\end{Lemma}

\subsection{Proof of Theorem 1}\label{s2.6}

Let $C_1$, $C_2$, $C_3$ and $C_4$ be disjoint oriented smooth closed
curves in $\R P^3$ generic in the sense of Theorem 1 (see Section
\ref{s1}). To prove Theorem 1, we have to find an orientation $\Gw$ of
0-manifold $\mathcal P(C_1,C_2,C_3,C_4)$ such that 
$$\sign(\mathcal P(C_1,C_2,C_3,C_4),\Gw)=
2\det \begin{pmatrix} \lk(C_1,C_2)&\lk(C_1,C_4)\\
                       \lk(C_3,C_2)&\lk(C_3,C_4)\end{pmatrix}.   $$

The genericity condition implies that any 
$L\in\mathcal P(C_1,C_2,C_3,C_4)$ can be turned into a pointed 
common secant of $C_1$, $C_2$ and $C_3$ in a unique way, 
and the pointed common secant is regular. 

Therefore by a small perturbation of $C_1$, $C_2$ and $C_3$ which is
trivial on a neighborhood of each common secant of all four curves, 
we can make the triple $C_1$, $C_2$, $C_3$ regular without changing 
$\mathcal P(C_1,C_2,C_3,C_4)$. Thus we can assume that triple of 
curves $C_1$, $C_2$, $C_3$ is regular from the very beginning. 

There is obvious bijection between $\mathcal P(C_1,C_2,C_3,C_4)$ and 
$C_4\cap U(\Sigma_{3,1})$. The intersection of $C_4$ with
$U(\Sigma_{3,1})$ is transversal by the fourth condition of general
position, see Section \ref{s1.1}. The orientations of $C_4$,
$\Sigma_{1,3}$ and the ambient space $\R P^3$ define an orientation of
the intersection $C_4\cap U(\Sigma_{3,1})$ turning it into an oriented
0-manifold. The intersection number $U(\Sigma_{1,3})\circ C_4$ 
is the signature of $C_4\cap U(\Sigma_{3,1})$.\qed

\subsection{In Affine Space}\label{s2.7}
Let $C_1$, $C_2$, $C_3$ be three smooth disjoint curves in $\R^3$.
Denote by $\mathcal A(C_1,C_2,C_3)$ the set of affine lines in $\R^3$
each of which passes through points $a\in C_1$, $b\in C_2$ and $c\in
C_3$ in this order. Put 
$ sa(C_1,C_2,C_3)=\bigcup_{L\in \mathcal A(C_,C_2,C_3)}L$.
Denote by $SA(C_1,C_2,C_3)$ the set
$$\{(x,p_1,p_2,p_3)\in\R^3\times C_1\times C_3\times
C_3\mid x,p_1,p_2,p_3 \text{ are colinear, }p_2\in[p_1,p_3]\}$$
and by 
$TA(C_1,C_2,C_3)$ the set $$\{(p_1,p_2,p_3)\in C_1\times C_2\times C_3 \mid  
 p_1,p_2,p_3 \text{ are colinear, }p_2\in[p_1,p_3]\}.$$ 

There are natural maps  
$$\begin{aligned}
ua:&\  SA(C_1,C_2,C_3)\to sa(C_1,C_2,C_3) : (x,p_1,p_2,p_3)\mapsto x\\
UA:&\  SA(C_1,C_2,C_3)\to \R^3 : (x,p_1,p_2,p_3)\mapsto x\\  
qa:&\  SA(C_1,C_2,C_3)\to TA(C_1,C_2,C_3) : (x,p_1,p_2,p_3)\mapsto (p_1,p_2,p_3)
\end{aligned}$$

Map $qa$ is a trivial line fibration  with trivialization 
$$SA(C_1,C_2,C_3)\to\R^1\times T(C_1,C_2,C_3),$$
under which the sections $P_1$,  $P_3$ defined by
 $$P_i=\{(x,p_1,p_2,p_3)\in S(C_1,C_2,C_3)\mid x=p_i\}.$$ 
are mapped to fibers 
$0\times T(C_1,C_2,C_3)$ and  $1\times T(C_1,C_2,C_3)$,
respectively.               

Inclusion $in:\R^3\subset\R P^3$ induces embeddings
$$\begin{aligned}
sa(C_1,C_2,C_3)&\to s(in(C_1),in(C_2),in(C_3))\\
SA(C_1,C_2,C_3)&\to S(in(C_1),in(C_2),in(C_3))\\
TA(C_1,C_2,C_3)&\to T(in(C_1),in(C_2),in(C_3))
\end{aligned}
$$
However, they are not necessarily surjective, because in the definition
of the spaces for affine situation there is condition
$p_2\in[p_1,p_3]$, which does not make sense and has no counter-part in 
the projective situation. Genericity assumptions considered in Section
\ref{s2.2} above, can be borrowed entirely (although one could reduce
them).
The orientation of $S(C_1,C_2,C_3)$ introduced in Section \ref{s2.3}
under the genericity assumptions can be repeated without changes in the
affine situation, or can be borrowed using the embedding 
$SA(C_1,C_2,C_3)\to S(in(C_1),in(C_2),in(C_3))$. 

Denote by $\Sigma$ the union of rays $[p_3,\infty)$ on the fibers of 
fibration $$qa:SA(C_1,C_2,C_3)\to TA(C_1,C_2,C_3).$$ 
Clearly, $\p\Sigma=P_3$. 

The following statement is a counter-part of Lemma \ref{L6}.

\begin{Lemma}\label{L9}
The degree of map $UA_3:P_3\to C_3$ defined by the natural map
$UA:SA(C_1,C_2,C_3)\to\R^3$ is equal to $\lk(C_1,C_2)$.
\end{Lemma} 

Proof is similar to the proof of Lemma \ref{L6} given above. We consider
the only point where the proofs differ.

The factor 2 disappeared, because the preimage of $p_3$ under $UA_3$
consists of points $(p_1,p_2,p_3,p_3)$ of pointed secants $(p_1,p_2,p_3)$
with $p_1\in C_1$, $p_2\in C_2$ and $p_2\in[p_1,p_3]$. Thus we count,
with appropriate signs, crossing points, where $C_2$ is above $C_1$, 
of the diagram for $C_1\cup C_2$ generated by
projection centered at $p_3$ . This gives the
linking number. In the projective case (in Lemma \ref{L6}) the condition 
 $p_2\in[p_1,p_3]$ was absent which resulted the factor 2.\qed

\begin{Lemma}\label{L10}
The image of the fundamental cycle of $\Sigma$ under mapping 
$UA:SA(C_1,C_2,C_3)\to\R^3$ is a smooth singular chain with closed
(non-compact) support and boundary $\lk(C_1,C_2)[C_3]$. 
\end{Lemma}

This is the counter-part of Lemma \ref{L7} and it follows immediately
from Lemma \ref{L9}.\qed

\begin{Lemma}\label{L11}
The intersection number of $UA_*[\Sigma]$ with an oriented closed
curve $C_4$ disjoint from $C_1\cup C_2\cup C_3$ is
$\lk(C_1,C_2)\lk(C_3,C_4)$.
\end{Lemma} \qed

After this point the proof of Theorem 2 runs like its counter-part
in Section \ref{s2.6}. The 0-dimensional manifold $\mathcal
A(C_1,C_2,C_3,C_4)$ is identified with $C_4\cap UA(\Sigma)$.

Theorems 3 and 4 are proved in the exactly same way, 
although they cannot be deduced formally from Theorem 2.

\subsection{Circles Meeting Curves}\label{s2.8}

\begin{Lemma}\label{L12}
For any disjoint oriented smooth closed curves $C_1$, $C_2$, $C_3$ and
$C_4$ in $\R^3$ and a point $p\in\R^3$ the set $\mathcal
S(C_1,C_2,C_3,C_4,p)$ of circles meeting $C_1$, $C_2$, $C_3$, $C_4$ and
$p$ in this cyclic order is a 0-manifold, which admits orientation $\Gw$
with signature $\lk(C_1,C_2)\lk(C_3,C_4)$.
\end{Lemma}

\begin{proof}Apply an inversion $I$ of $R^3$ centered at $p$. This turns the
circles belonging to  $\mathcal S(C_1,C_2,C_3,C_4,p)$ to lines belonging
to $\mathcal A(I(C_1),I(C_2),I(C_3), I(C_4))$. Now use Theorem 2.
\end{proof}

For a generic collection $C_1$, $C_2$, $C_3$, $C_4$ and $C_5$ of
smooth closed curves in $\R^3$ denote by
$SC(C_1,C_2,C_3,C_4,C_5)$ the set of $(x,p_1,p_2,p_3,p_4,p_5)\in
\R^3\times C_1\times C_2\times C_3\times C_4\times C_5$ such that there 
exists a circle passing through 
$p_1,p_2,p_3,p_4,p_5$  in this cyclic order and $x$ lies on this circle,
too.

Denote by $TC(C_1,C_2,C_3,C_4,C_5)$ the set of $(p_1,p_2,p_3,p_4,p_5)\in
C_1\times C_2\times C_3\times C_4\times C_5$ such that there 
exists a circle passing through $p_1,p_2,p_3,p_4,p_5$  in this cyclic
order. 

Formula $(x,p_1,p_2,p_3,p_4,p_5)\mapsto (p_1,p_2,p_3,p_4,p_5)$ defines a
fibration 
$$SC(C_1,C_2,C_3,C_4,C_5)\to TC(C_1,C_2,C_3,C_4,C_5)$$ 
with fiber circle. It has five disjoint sections defined by formulas
$$(p_1,p_2,p_3,p_4,p_5)\mapsto (p_i,p_1,p_2,p_3,p_4,p_5)\text{ with }
i=1,\dots,5.$$ 
The images of these sections can be described by the following formulas:
$$P_i=\{(x,p_1,p_2,p_3,p_4,p_5)\in SC(C_1,C_2,C_3,C_4,C_5)\mid x=p_i\}. $$
The part of $SC(C_1,C_2,C_3,C_4,C_5)$ bounded by sections $P_i$ and
$P_{i+1}$ and disjoint with other $P_j$ is denoted by $\Sigma_{i}$
(here, when we write $i+1$, we treat $i$ modulo 5).

Similarly to what was done in Section \ref{s2.3}, for generic
collection of curves the space $SC(C_1,C_2,C_3,C_4,C_5)$ is equipped
with a natural orientation. Equip sections $P_i$ with orientations such
that $\p\Sigma_i=P_i\cup(-P_{i+1})$.

\begin{Lemma}\label{L13}
The degree of map $UC_iP_i\to C_i$ defined by the natural map 
$$UC:SC(C_1,C_2,C_3,C_4,C_5)\to\R^3:(x,p_1,p_2,p_3,p_4,p_5)\mapsto x$$
is equal to $\lk(C_{i+1},C_{i+2})\lk(C_{i+3},C_{i+4})$ (here again the
indices are treated as elements of $\Z/5\Z$).
\end{Lemma}

This lemma is similar to Lemmas \ref{L6} and \ref{L9}.
The proof of this lemma is based on Lemma \ref{L12}.

After this point our proof of Theorem 5 goes along the same scheme as
proofs of Theorems 1 and 2 above. We identify $\mathcal
S(C_1,C_2,C_3,C_4,C_5,C_6$ with $C_6\cap UC(\Sigma_5)$ and orient the
latter using the given orientation of $C_6$ and the orientation of
$\Sigma_5$.

Theorems 6 and 7 are proved by the same arguments, but applied to spaces
appropriately changed.

\end{document}

%% file: 4lines+.pstex_t
\begin{picture}(0,0)%
\includegraphics{4lines+.pstex}%
\end{picture}%
\setlength{\unitlength}{4144sp}%
\begingroup\makeatletter\ifx\SetFigFont\undefined%
\gdef\SetFigFont#1#2#3#4#5{%
  \reset@font\fontsize{#1}{#2pt}%
  \fontfamily{#3}\fontseries{#4}\fontshape{#5}%
  \selectfont}%
\fi\endgroup%
\begin{picture}(1585,1939)(76,-1120)
\put( 91,-1051){\makebox(0,0)[lb]{\smash{{\SetFigFont{12}{14.4}{\familydefault}{\mddefault}{\updefault}{\color[rgb]{0,0,0}$C_1$}%
}}}}
\put(496,-1051){\makebox(0,0)[lb]{\smash{{\SetFigFont{12}{14.4}{\familydefault}{\mddefault}{\updefault}{\color[rgb]{0,0,0}$C_2$}%
}}}}
\put(1036,-1051){\makebox(0,0)[lb]{\smash{{\SetFigFont{12}{14.4}{\familydefault}{\mddefault}{\updefault}{\color[rgb]{0,0,0}$C_3$}%
}}}}
\put(1531,-1051){\makebox(0,0)[lb]{\smash{{\SetFigFont{12}{14.4}{\familydefault}{\mddefault}{\updefault}{\color[rgb]{0,0,0}$C_4$}%
}}}}
\end{picture}%

%% file: 4lines-.pstex_t
\begin{picture}(0,0)%
\includegraphics{4lines-.pstex}%
\end{picture}%
\setlength{\unitlength}{4144sp}%
\begingroup\makeatletter\ifx\SetFigFont\undefined%
\gdef\SetFigFont#1#2#3#4#5{%
  \reset@font\fontsize{#1}{#2pt}%
  \fontfamily{#3}\fontseries{#4}\fontshape{#5}%
  \selectfont}%
\fi\endgroup%
\begin{picture}(1540,1954)(121,-1120)
\put(136,-1051){\makebox(0,0)[lb]{\smash{{\SetFigFont{12}{14.4}{\familydefault}{\mddefault}{\updefault}{\color[rgb]{0,0,0}$C_1$}%
}}}}
\put(541,-1051){\makebox(0,0)[lb]{\smash{{\SetFigFont{12}{14.4}{\familydefault}{\mddefault}{\updefault}{\color[rgb]{0,0,0}$C_2$}%
}}}}
\put(1036,-1051){\makebox(0,0)[lb]{\smash{{\SetFigFont{12}{14.4}{\familydefault}{\mddefault}{\updefault}{\color[rgb]{0,0,0}$C_3$}%
}}}}
\put(1486,-1051){\makebox(0,0)[lb]{\smash{{\SetFigFont{12}{14.4}{\familydefault}{\mddefault}{\updefault}{\color[rgb]{0,0,0}$C_4$}%
}}}}
\end{picture}%

%% file: 4lines0.pstex_t
\begin{picture}(0,0)%
\includegraphics{4lines0.pstex}%
\end{picture}%
\setlength{\unitlength}{4144sp}%
\begingroup\makeatletter\ifx\SetFigFont\undefined%
\gdef\SetFigFont#1#2#3#4#5{%
  \reset@font\fontsize{#1}{#2pt}%
  \fontfamily{#3}\fontseries{#4}\fontshape{#5}%
  \selectfont}%
\fi\endgroup%
\begin{picture}(1585,1939)(76,-1120)
\put( 91,-1051){\makebox(0,0)[lb]{\smash{{\SetFigFont{12}{14.4}{\familydefault}{\mddefault}{\updefault}{\color[rgb]{0,0,0}$C_1$}%
}}}}
\put(541,-1051){\makebox(0,0)[lb]{\smash{{\SetFigFont{12}{14.4}{\familydefault}{\mddefault}{\updefault}{\color[rgb]{0,0,0}$C_2$}%
}}}}
\put(1081,-1051){\makebox(0,0)[lb]{\smash{{\SetFigFont{12}{14.4}{\familydefault}{\mddefault}{\updefault}{\color[rgb]{0,0,0}$C_3$}%
}}}}
\put(1576,-1051){\makebox(0,0)[lb]{\smash{{\SetFigFont{12}{14.4}{\familydefault}{\mddefault}{\updefault}{\color[rgb]{0,0,0}$C_4$}%
}}}}
\end{picture}%

%% file: 4linesOnH+.pstex_t
\begin{picture}(0,0)%
\includegraphics{4linesOnH+.pstex}%
\end{picture}%
\setlength{\unitlength}{4144sp}%
\begingroup\makeatletter\ifx\SetFigFont\undefined%
\gdef\SetFigFont#1#2#3#4#5{%
  \reset@font\fontsize{#1}{#2pt}%
  \fontfamily{#3}\fontseries{#4}\fontshape{#5}%
  \selectfont}%
\fi\endgroup%
\begin{picture}(5455,2489)(211,-1615)
\put(3736,-1546){\makebox(0,0)[lb]{\smash{{\SetFigFont{12}{14.4}{\familydefault}{\mddefault}{\updefault}{\color[rgb]{0,0,0}$C_2$}%
}}}}
\put(4546,-1546){\makebox(0,0)[lb]{\smash{{\SetFigFont{12}{14.4}{\familydefault}{\mddefault}{\updefault}{\color[rgb]{0,0,0}$C_3$}%
}}}}
\put(5221,-1546){\makebox(0,0)[lb]{\smash{{\SetFigFont{12}{14.4}{\familydefault}{\mddefault}{\updefault}{\color[rgb]{0,0,0}$C_4$}%
}}}}
\put(3196,-1546){\makebox(0,0)[lb]{\smash{{\SetFigFont{12}{14.4}{\familydefault}{\mddefault}{\updefault}{\color[rgb]{0,0,0}$C_1$}%
}}}}
\put(766,-1546){\makebox(0,0)[lb]{\smash{{\SetFigFont{12}{14.4}{\familydefault}{\mddefault}{\updefault}{\color[rgb]{0,0,0}$C_2$}%
}}}}
\put(1576,-1546){\makebox(0,0)[lb]{\smash{{\SetFigFont{12}{14.4}{\familydefault}{\mddefault}{\updefault}{\color[rgb]{0,0,0}$C_3$}%
}}}}
\put(2251,-1546){\makebox(0,0)[lb]{\smash{{\SetFigFont{12}{14.4}{\familydefault}{\mddefault}{\updefault}{\color[rgb]{0,0,0}$C_4$}%
}}}}
\put(226,-1546){\makebox(0,0)[lb]{\smash{{\SetFigFont{12}{14.4}{\familydefault}{\mddefault}{\updefault}{\color[rgb]{0,0,0}$C_1$}%
}}}}
\end{picture}%

%% file: lineOffH.pstex_t
\begin{picture}(0,0)%
\includegraphics{lineOffH.pstex}%
\end{picture}%
\setlength{\unitlength}{3947sp}%
\begingroup\makeatletter\ifx\SetFigFont\undefined%
\gdef\SetFigFont#1#2#3#4#5{%
  \reset@font\fontsize{#1}{#2pt}%
  \fontfamily{#3}\fontseries{#4}\fontshape{#5}%
  \selectfont}%
\fi\endgroup%
\begin{picture}(5794,4361)(61,-3580)
\put(1426,-3511){\makebox(0,0)[lb]{\smash{{\SetFigFont{12}{14.4}{\familydefault}{\mddefault}{\updefault}{\color[rgb]{0,0,0}(d)}%
}}}}
\put(3676,-3511){\makebox(0,0)[lb]{\smash{{\SetFigFont{12}{14.4}{\familydefault}{\mddefault}{\updefault}{\color[rgb]{0,0,0}(e)}%
}}}}
\put(4426,-3286){\makebox(0,0)[lb]{\smash{{\SetFigFont{12}{14.4}{\familydefault}{\mddefault}{\updefault}{\color[rgb]{0,0,0}$H$}%
}}}}
\put(2101,-3286){\makebox(0,0)[lb]{\smash{{\SetFigFont{12}{14.4}{\familydefault}{\mddefault}{\updefault}{\color[rgb]{0,0,0}$H$}%
}}}}
\put(2026,-2011){\makebox(0,0)[lb]{\smash{{\SetFigFont{12}{14.4}{\familydefault}{\mddefault}{\updefault}{\color[rgb]{0,0,0}$M$}%
}}}}
\put(4276,-2011){\makebox(0,0)[lb]{\smash{{\SetFigFont{12}{14.4}{\familydefault}{\mddefault}{\updefault}{\color[rgb]{0,0,0}$M$}%
}}}}
\put(4501,-1411){\makebox(0,0)[lb]{\smash{{\SetFigFont{12}{14.4}{\familydefault}{\mddefault}{\updefault}{\color[rgb]{0,0,0}$N$}%
}}}}
\put(901,-1711){\makebox(0,0)[lb]{\smash{{\SetFigFont{12}{14.4}{\familydefault}{\mddefault}{\updefault}{\color[rgb]{0,0,0}$Q$}%
}}}}
\put(3151,-1711){\makebox(0,0)[lb]{\smash{{\SetFigFont{12}{14.4}{\familydefault}{\mddefault}{\updefault}{\color[rgb]{0,0,0}$Q$}%
}}}}
\put(3226,-2236){\makebox(0,0)[lb]{\smash{{\SetFigFont{12}{14.4}{\familydefault}{\mddefault}{\updefault}{\color[rgb]{0,0,0}$p$}%
}}}}
\put(526,-1111){\makebox(0,0)[lb]{\smash{{\SetFigFont{12}{14.4}{\familydefault}{\mddefault}{\updefault}{\color[rgb]{0,0,0}(a)}%
}}}}
\put(2776,-1111){\makebox(0,0)[lb]{\smash{{\SetFigFont{12}{14.4}{\familydefault}{\mddefault}{\updefault}{\color[rgb]{0,0,0}(b)}%
}}}}
\put(4876,-1111){\makebox(0,0)[lb]{\smash{{\SetFigFont{12}{14.4}{\familydefault}{\mddefault}{\updefault}{\color[rgb]{0,0,0}(c)}%
}}}}
\put(1201,-961){\makebox(0,0)[lb]{\smash{{\SetFigFont{12}{14.4}{\familydefault}{\mddefault}{\updefault}{\color[rgb]{0,0,0}$H$}%
}}}}
\put(3376,-961){\makebox(0,0)[lb]{\smash{{\SetFigFont{12}{14.4}{\familydefault}{\mddefault}{\updefault}{\color[rgb]{0,0,0}$H$}%
}}}}
\put(5476,-961){\makebox(0,0)[lb]{\smash{{\SetFigFont{12}{14.4}{\familydefault}{\mddefault}{\updefault}{\color[rgb]{0,0,0}$H$}%
}}}}
\put(526,-61){\makebox(0,0)[lb]{\smash{{\SetFigFont{12}{14.4}{\familydefault}{\mddefault}{\updefault}{\color[rgb]{0,0,0}$L$}%
}}}}
\put(2701,-61){\makebox(0,0)[lb]{\smash{{\SetFigFont{12}{14.4}{\familydefault}{\mddefault}{\updefault}{\color[rgb]{0,0,0}$L$}%
}}}}
\put(4801,-61){\makebox(0,0)[lb]{\smash{{\SetFigFont{12}{14.4}{\familydefault}{\mddefault}{\updefault}{\color[rgb]{0,0,0}$L$}%
}}}}
\put(5401,314){\makebox(0,0)[lb]{\smash{{\SetFigFont{12}{14.4}{\familydefault}{\mddefault}{\updefault}{\color[rgb]{0,0,0}$M$}%
}}}}
\put(2101,-211){\makebox(0,0)[lb]{\smash{{\SetFigFont{12}{14.4}{\familydefault}{\mddefault}{\updefault}{\color[rgb]{0,0,0}$P$}%
}}}}
\put(4201,-211){\makebox(0,0)[lb]{\smash{{\SetFigFont{12}{14.4}{\familydefault}{\mddefault}{\updefault}{\color[rgb]{0,0,0}$P$}%
}}}}
\end{picture}%

%% file: 4lines0andH.pstex_t
\begin{picture}(0,0)%
\includegraphics{4lines0andH.pstex}%
\end{picture}%
\setlength{\unitlength}{4144sp}%
\begingroup\makeatletter\ifx\SetFigFont\undefined%
\gdef\SetFigFont#1#2#3#4#5{%
  \reset@font\fontsize{#1}{#2pt}%
  \fontfamily{#3}\fontseries{#4}\fontshape{#5}%
  \selectfont}%
\fi\endgroup%
\begin{picture}(3305,2487)(-31,-1615)
\put(766,-1546){\makebox(0,0)[lb]{\smash{{\SetFigFont{12}{14.4}{\familydefault}{\mddefault}{\updefault}{\color[rgb]{0,0,0}$C_2$}%
}}}}
\put(1576,-1546){\makebox(0,0)[lb]{\smash{{\SetFigFont{12}{14.4}{\familydefault}{\mddefault}{\updefault}{\color[rgb]{0,0,0}$C_3$}%
}}}}
\put(226,-1546){\makebox(0,0)[lb]{\smash{{\SetFigFont{12}{14.4}{\familydefault}{\mddefault}{\updefault}{\color[rgb]{0,0,0}$C_1$}%
}}}}
\put(3061,-1546){\makebox(0,0)[lb]{\smash{{\SetFigFont{12}{14.4}{\familydefault}{\mddefault}{\updefault}{\color[rgb]{0,0,0}$C_4$}%
}}}}
\end{picture}%

%% file: 3+1curves.pstex_t
\begin{picture}(0,0)%
\includegraphics{3+1curves.pstex}%
\end{picture}%
\setlength{\unitlength}{4144sp}%
\begingroup\makeatletter\ifx\SetFigFont\undefined%
\gdef\SetFigFont#1#2#3#4#5{%
  \reset@font\fontsize{#1}{#2pt}%
  \fontfamily{#3}\fontseries{#4}\fontshape{#5}%
  \selectfont}%
\fi\endgroup%
\begin{picture}(2332,2087)(34,-1215)
\put(1546,-1041){\makebox(0,0)[lb]{\smash{{\SetFigFont{12}{14.4}{\rmdefault}{\mddefault}{\updefault}{\color[rgb]{0,0,.82}$C_4$}%
}}}}
\put(1006,-1146){\makebox(0,0)[lb]{\smash{{\SetFigFont{12}{14.4}{\rmdefault}{\mddefault}{\updefault}{\color[rgb]{0,0,.82}$C_1$}%
}}}}
\put(1736,-351){\makebox(0,0)[lb]{\smash{{\SetFigFont{12}{14.4}{\rmdefault}{\mddefault}{\updefault}{\color[rgb]{0,0,.82}$C_2$}%
}}}}
\put(2351,549){\makebox(0,0)[lb]{\smash{{\SetFigFont{12}{14.4}{\rmdefault}{\mddefault}{\updefault}{\color[rgb]{0,0,0}$C_3$}%
}}}}
\put(1146,-111){\makebox(0,0)[lb]{\smash{{\SetFigFont{12}{14.4}{\rmdefault}{\mddefault}{\updefault}{\color[rgb]{0,0,0}$x$}%
}}}}
\put(2121,154){\makebox(0,0)[lb]{\smash{{\SetFigFont{12}{14.4}{\rmdefault}{\mddefault}{\updefault}{\color[rgb]{0,0,.82}$s(C_1,C_2,C_3)$}%
}}}}
\end{picture}%

%% file: reg_sec.pstex_t
\begin{picture}(0,0)%
\includegraphics{reg_sec.pstex}%
\end{picture}%
\setlength{\unitlength}{3315sp}%
\begingroup\makeatletter\ifx\SetFigFont\undefined%
\gdef\SetFigFont#1#2#3#4#5{%
  \reset@font\fontsize{#1}{#2pt}%
  \fontfamily{#3}\fontseries{#4}\fontshape{#5}%
  \selectfont}%
\fi\endgroup%
\begin{picture}(2278,2384)(431,-1448)
\put(936,594){\makebox(0,0)[lb]{\smash{{\SetFigFont{10}{12.0}{\familydefault}{\mddefault}{\updefault}{\color[rgb]{0,0,0}$L_1$}%
}}}}
\put(1551,597){\makebox(0,0)[lb]{\smash{{\SetFigFont{10}{12.0}{\familydefault}{\mddefault}{\updefault}{\color[rgb]{0,0,0}$L_2$}%
}}}}
\put(2202,596){\makebox(0,0)[lb]{\smash{{\SetFigFont{10}{12.0}{\familydefault}{\mddefault}{\updefault}{\color[rgb]{0,0,0}$L_3$}%
}}}}
\end{picture}%

%% file: alm_reg_sec.pstex_t
\begin{picture}(0,0)%
\includegraphics{alm_reg_sec.pstex}%
\end{picture}%
\setlength{\unitlength}{3108sp}%
\begingroup\makeatletter\ifx\SetFigFont\undefined%
\gdef\SetFigFont#1#2#3#4#5{%
  \reset@font\fontsize{#1}{#2pt}%
  \fontfamily{#3}\fontseries{#4}\fontshape{#5}%
  \selectfont}%
\fi\endgroup%
\begin{picture}(5063,2873)(58,-1755)
\put(3041,-1101){\makebox(0,0)[lb]{\smash{{\SetFigFont{9}{10.8}{\rmdefault}{\mddefault}{\updefault}{\color[rgb]{0,.56,0}view from $p_k$ along}%
}}}}
\put(3041,-1326){\makebox(0,0)[lb]{\smash{{\SetFigFont{9}{10.8}{\rmdefault}{\mddefault}{\updefault}{\color[rgb]{0,.56,0}almost regular secant}%
}}}}
\put(2736,669){\makebox(0,0)[lb]{\smash{{\SetFigFont{9}{10.8}{\rmdefault}{\mddefault}{\updefault}{\color[rgb]{.56,0,.56}plane containing $L_i$ and $L_j$}%
}}}}
\put(5106,-506){\makebox(0,0)[lb]{\smash{{\SetFigFont{9}{10.8}{\rmdefault}{\mddefault}{\updefault}{\color[rgb]{0,0,0}$L_i=L_k$}%
}}}}
\put(4796, 64){\makebox(0,0)[lb]{\smash{{\SetFigFont{9}{10.8}{\rmdefault}{\mddefault}{\updefault}{\color[rgb]{0,0,.82}$C_i$}%
}}}}
\put(4821,-911){\makebox(0,0)[lb]{\smash{{\SetFigFont{9}{10.8}{\rmdefault}{\mddefault}{\updefault}{\color[rgb]{0,0,.82}$C_j$}%
}}}}
\put(1271,959){\makebox(0,0)[lb]{\smash{{\SetFigFont{9}{10.8}{\rmdefault}{\mddefault}{\updefault}{\color[rgb]{0,.56,0}almost regular secant}%
}}}}
\put(1191, 39){\makebox(0,0)[lb]{\smash{{\SetFigFont{9}{10.8}{\rmdefault}{\mddefault}{\updefault}{\color[rgb]{0,.56,0}$p_j$}%
}}}}
\put(2136,459){\makebox(0,0)[lb]{\smash{{\SetFigFont{9}{10.8}{\rmdefault}{\mddefault}{\updefault}{\color[rgb]{0,.56,0}$L_j$}%
}}}}
\put(2021,-531){\makebox(0,0)[lb]{\smash{{\SetFigFont{9}{10.8}{\rmdefault}{\mddefault}{\updefault}{\color[rgb]{0,.56,0}$L_i$}%
}}}}
\put(511,329){\makebox(0,0)[lb]{\smash{{\SetFigFont{9}{10.8}{\rmdefault}{\mddefault}{\updefault}{\color[rgb]{0,.56,0}$C_j$}%
}}}}
\put(1281,-291){\makebox(0,0)[lb]{\smash{{\SetFigFont{9}{10.8}{\rmdefault}{\mddefault}{\updefault}{\color[rgb]{0,.56,0}$p_i$}%
}}}}
\put(2036,-996){\makebox(0,0)[lb]{\smash{{\SetFigFont{9}{10.8}{\rmdefault}{\mddefault}{\updefault}{\color[rgb]{0,.56,0}$C_i$}%
}}}}
\put(511,-1066){\makebox(0,0)[lb]{\smash{{\SetFigFont{9}{10.8}{\rmdefault}{\mddefault}{\updefault}{\color[rgb]{0,0,.82}$L_k$}%
}}}}
\put(1446,-1266){\makebox(0,0)[lb]{\smash{{\SetFigFont{9}{10.8}{\rmdefault}{\mddefault}{\updefault}{\color[rgb]{0,.56,0}$p_k$}%
}}}}
\put(226,-1351){\makebox(0,0)[lb]{\smash{{\SetFigFont{9}{10.8}{\rmdefault}{\mddefault}{\updefault}{\color[rgb]{0,0,.82}$C_k$}%
}}}}
\end{picture}%

%% file: whitn_umbr.pstex_t
\begin{picture}(0,0)%
\includegraphics{whitn_umbr.pstex}%
\end{picture}%
\setlength{\unitlength}{3315sp}%
\begingroup\makeatletter\ifx\SetFigFont\undefined%
\gdef\SetFigFont#1#2#3#4#5{%
  \reset@font\fontsize{#1}{#2pt}%
  \fontfamily{#3}\fontseries{#4}\fontshape{#5}%
  \selectfont}%
\fi\endgroup%
\begin{picture}(2936,2052)(-199,-1180)
\put(1421,639){\makebox(0,0)[lb]{\smash{{\SetFigFont{10}{12.0}{\rmdefault}{\mddefault}{\updefault}{\color[rgb]{0,.56,.56}$C_1$}%
}}}}
\put(-184,-326){\makebox(0,0)[lb]{\smash{{\SetFigFont{10}{12.0}{\rmdefault}{\mddefault}{\updefault}{\color[rgb]{0,.56,.56}$C_2$}%
}}}}
\put(2631,-291){\makebox(0,0)[lb]{\smash{{\SetFigFont{10}{12.0}{\rmdefault}{\mddefault}{\updefault}{\color[rgb]{0,.56,.56}$C_3$}%
}}}}
\put(381, 99){\makebox(0,0)[lb]{\smash{{\SetFigFont{10}{12.0}{\rmdefault}{\mddefault}{\updefault}{\color[rgb]{0,.56,.56}almost regular secant}%
}}}}
\put(1046,379){\makebox(0,0)[lb]{\smash{{\SetFigFont{10}{12.0}{\rmdefault}{\mddefault}{\updefault}{\color[rgb]{0,.56,.56}$p_1$}%
}}}}
\put( 56,389){\makebox(0,0)[lb]{\smash{{\SetFigFont{10}{12.0}{\rmdefault}{\mddefault}{\updefault}{\color[rgb]{0,.56,.56}$p_2$}%
}}}}
\put(2141,384){\makebox(0,0)[lb]{\smash{{\SetFigFont{10}{12.0}{\rmdefault}{\mddefault}{\updefault}{\color[rgb]{0,.56,.56}$p_3$}%
}}}}
\end{picture}%